\documentclass[letterpaper, 10 pt, conference]{ieeeconf}%

\usepackage[cmex10]{amsmath}

\usepackage{ifpdf}

\usepackage{algorithm}
\usepackage{algpseudocode}
\usepackage{color}
\usepackage{url}
\ifpdf
\usepackage[pdftex]{graphicx}
\else
\usepackage[dvips]{graphicx}
\fi
\usepackage{tabularx}
\DeclareGraphicsExtensions{.eps}
\usepackage[caption=false,font=footnotesize]{subfig}
\usepackage{comment}
\usepackage{bm}
\usepackage{amsmath,amsthm,amsfonts,amssymb}
\usepackage[noadjust]{cite}
\usepackage{nomencl}
\theoremstyle{definition}
\newtheorem{theorem}{\textbf{Theorem}}

\newtheorem{proposition}[theorem]{\textbf{Proposition}}
\newtheorem{remark}[theorem]{\textbf{Remark}}

\newcommand{\up}[1]{^\mathrm{#1}}
\IEEEoverridecommandlockouts
\begin{document}
%
% paper title
% Titles are generally capitalized except for words such as a, an, and, as,
% at, but, by, for, in, nor, of, on, or, the, to and up, which are usually
% not capitalized unless they are the first or last word of the title.
% Linebreaks \\ can be used within to get better formatting as desired.
% Do not put math or special symbols in the title.
\title{Operational Valuation of Energy Storage under Multi-stage Price Uncertainties}

%% To specify the authors when (number of affiliations <= 2)
\author{Bolun~Xu,~\IEEEmembership{Member,~IEEE,}
	% \vspace{-.7cm}
	Magnus Korp\aa s,~\IEEEmembership{Member,~IEEE,}
	% \vspace{-.7cm}
	Audun Botterud,~\IEEEmembership{Member,~IEEE,}
	\thanks{B.~Xu is with Columbia University, NY, USA. M. Korp\aa s is with Norwegian University of Science and Technology, Trondheim, Norway. A.~Botterud is with MIT, MA, USA. }
}

%% To specify the authors when (number of affiliations > 2)
% \author{\IEEEauthorblockN{Author n.1\IEEEauthorrefmark{1},
% Author n.2\IEEEauthorrefmark{2},
% Author n.3\IEEEauthorrefmark{3}, 
% Author n.4\IEEEauthorrefmark{3} and
% Author n.5\IEEEauthorrefmark{4}}
% \IEEEauthorblockA{\IEEEauthorrefmark{1} Department Name of Organization A\\
% Name of the organization A,
% Address A\\ Emails if wanted}
% \IEEEauthorblockA{\IEEEauthorrefmark{2} Department Name of Organization B\\
% Name of the organization B,
% Address B\\ Emails if wanted}
% \IEEEauthorblockA{\IEEEauthorrefmark{3} Department Name of Organization C\\
% Name of the organization C,
% Address C\\ Emails if wanted}
% \IEEEauthorblockA{\IEEEauthorrefmark{4}Department Name of Organization D\\
% Name of the organization D,
% Address D\\ Emails if wanted}
% }

% make the title area
\maketitle

% As a general rule, do not put math, special symbols or citations
% in the abstract
\begin{abstract}
This paper presents an analytical method for calculating the operational value of an energy storage device under multi-stage price uncertainties. Our solution calculates the storage value function from price distribution functions directly instead of sampling discrete scenarios, offering improved modeling accuracy over tail distribution events such as price spikes and negative prices. The analytical algorithm offers very high computational efficiency in solving multi-stage stochastic programming for energy storage and can easily be implemented within any software and hardware platform, while numerical simulation results show the proposed method is up to 100,000 times faster than a benchmark stochastic-dual dynamic programming solver even in small test cases. Case studies are included to demonstrate the impact of price variability on the valuation results, and a battery charging example using historical prices for New York City. 
\end{abstract}


\IEEEpeerreviewmaketitle
% Use this to place sponsorships
% \thanksto{Applicable sponsors, if any, should be placed using the \emph{thanksto} command.}

\section{Introduction}
Electricity markets are lowering participating barriers for energy storage, and many system operators have proposed new market policies for storage participants to bid according to their own economic valuation~\cite{sakti2018review}.  However, the operating value of storage devices depends on both the current and future system conditions due to their limited energy capacity, making their valuation substantially different from conventional thermal generators which are primarily based on fuel costs~\cite{kirschen2018fundamentals}. 
In addition, the valuation must also account for the system variability and the occurrence of sudden deviations that often results in price spikes~\cite{kirschen2003demand} or negative prices~\cite{martinez2016impact}. Operational planning for conventional generators are less focused on variability as generators are unlikely to be constrained on fuel storage and are often unable to response to these extreme prices due to their constrained ramp rate and start-up limits. In contrast, batteries can ramp from idle to full capacity within a few seconds, thus it is crucial for batteries to position their storage level for responding to sudden system imbalances according to the price signal, maximizing their operating profits while contributing to system security. 

The most convenient way of valuating storage is through price forecasts, as all detailed system information is oftentimes not available to market participants due to confidentiality. 
Electricity prices are typically forecastd as stage-wise independent processes, e.g. using time series analysis such as auto-regressive integrated moving average (ARIMA) ~\cite{shapiro2011analysis} for day-ahead price prediction~\cite{conejo2005day}, which calculates the expected price and the error distribution for each time period. The real-time price can be modeled as error functions on top of the day-ahead price published ahead of time by the system operator~\cite{tang2016model}. The occurrence of extreme price spikes can also be designed as features in time series forecast models or characterized as tail events in price error distributions~\cite{weron2008forecasting}. Recent studies~\cite{lakshminarayanan2017simple} also show promises in applying deep neural networks for probabilistic forecasts.

Energy storage arbitrage take advantage of price differences and it is therefore crucial to model the multi-stage price uncertainty in operational valuations, especially the distribution of tail events such as price spikes and negative prices. Optimizing energy storage operation under multi-stage uncertainties have been explored in varies applications using methods such as stochastic dynamic programming~\cite{xi2014stochastic}, stochastic dual dynamic programming~\cite{megel2015stochastic}, back propagation~\cite{salas2017benchmarking}, approximate dynamic programming~\cite{jiang2015optimal}, and reinforcement learning~\cite{wang2018energy}. These methods require discretization of the probability, action, and state spaces, which makes it difficult to model the impact of tail events in the distribution. We propose an analytical approach to solve the multi-stage price arbitrage problem for energy storage which obtains the storage value under instantaneous response to new price realizations. The main contributions of the paper are listed as follows:
\begin{enumerate}
    \item Our approach is based on price distribution functions directly instead of having to discretize distribution samples, offering better modeling accuracy for distribution tail events such as price spikes and negative prices. 
    \item Our approach calculates the storage value function analytically from the distribution function and the value function from the next period, providing very high computational efficiency. 
    \item Our approach solves the multi-stage energy storage arbitrage problem under linear time complexity and constant space complexity, offering almost instant computation over hundreds of forecast periods. 
    \item We provide case studies for energy storage operation using New York City prices.
\end{enumerate}
The rest of the paper is organized as follows: Section~II formulates the optimization problem. Section~III provides the theoretical results and the solution algorithm. Simulation results are listed in Section~IV and Section~V concludes the paper.

\section{Formulation}
\subsection{Formulation}
Our valuation framework takes a price-taker perspective by assuming that the storage operation will have no impact on the price formation. We model the electricity price $\lambda_t\in\mathbb{R}$ over time period $t$ as a stage-wise independent random process with probability distribution function $f_t(\cdot)$ and cumulative distribution function $F_t(\cdot)$. % Note that we assume each time period has a unique price distribution, thus time dependent information such as day-ahead prices (published before real-time prices) can be used to form the real-time price distribution~\cite{tang2016model}. 
We model the storage operation as a nonanticipatory control policy~\cite{shapiro2009lectures}, i.e., the storage control over a time interval will not be dependent on future price realizations. The objective of this control policy is to maximize the expected price arbitrage profit, which is equivalent to maximizing the social welfare under the price-taker setting~\cite{castillo2013profit}. The optimal nonanticipatory policy knowing the price distribution over a future period $\{\lambda_t\;|\;t\in\{1,\dotsc,T\}\}$ can be formulated as
\begin{subequations}\label{p1}
\begin{gather}
    \max_{p_t} \mathbb{E}\Big[\sum_{t=1}^T \lambda_t (p\up{d}_t-p\up{c}_t) - cp\up{d}_t \Big]  + V_T(e_T)\label{p1_obj}
\end{gather}
subjects to 
\begin{gather}
\{p\up{d}_t, p\up{c}_t\} \in \text{Nonanticipatory policies} \label{p1_c4} \\
    0 \leq p\up{c}_t \leq P,\; 0\leq p\up{d}_t \leq P \label{p1_c2} \\
    \text{$p\up{c}_t$ or $p\up{d}_t$ is zero for any $t$} \label{p1_c5}\\
    e_t - e_{t-1} = -p\up{d}_t/\eta + p\up{c}_t\eta \label{p1_c1}\\
    0 \leq e_t \leq E \label{p1_c3}
\end{gather}
\end{subequations}
where $p\up{c}_t$ is the energy charged into the storage over time period $t$, and $p\up{d}_t$ is the discharged energy over period $t$, \eqref{p1_c4} enforces the control to be nonanticipatory, and \eqref{p1_c5} enforces that storage cannot charge and discharge during the same time period. 
$e_t$ is the storage state-of-charge (SoC) over period $t$ modeled in \eqref{p1_c1} as the charge and discharge energy subject to the efficiency $\eta$, and \eqref{p1_c2} enforces the upper and lower energy bound. Note that we have normalized the time period duration into $p\up{c}_t$ and $p\up{d}_t$ so no duration coefficient is included here for presentation simplicity. The optimization maximizes the expected profit considering price and the marginal discharge cost $c$, which represents storage operation and management costs such as degradation~\cite{xu2017factoring}. At last, the final storage SoC is influenced by the end-value function $V_T(\cdot)$, which for instance can be utilized to ensure that electric vehicles are sufficiently charged at the end of the period (by designing $V_T$ as an indicator function based on the target final SoC level).

% \begin{enumerate}
%     \item $\Pi$ is the set of all feasible policies that satisfy constraints \eqref{p1_c1}--\eqref{p1_c3} and the nonanticipativity constraint in which the control result shall not depend on the future realization of the price.
%     \item $\bm{p}^{\pi} = \{p\up{\pi}_t\}$ is the battery operation controlled by the policy $\pi$.
%     \item $e_t$ is the state-of-charge.
%     \item $C(\cdot)$ is the cost function for operation the battery. 
%     \item $V_T(\cdot)$ is the value of the ending state-of-charge $e_t$ over the control period. For example, $V_T(\cdot)$ could be a step-function to enforce the battery being at least $80\%$ full at the end in the application of arbitraging a plug-in vehicle. 
% \end{enumerate}

% \section{Main Results}

\subsection{Value Function and Stochastic Dynamic Programming}
We use stochastic dynamic programming to solve the nonanticipatory constraint in problem \eqref{p1}, by working backwards and recursively solving a single-period optimization  ($\forall t<T$):
\begin{subequations}
\begin{align}
    Q_{t-1}(e_{t-1}) &= \max_{p\up{c}_t, p\up{d}_t} \lambda_t (p\up{d}_t-p\up{c}_t) - cp\up{d}_t + V_t(e_t) \nonumber\\
    & \text{ s.t. \eqref{p1_c2}--\eqref{p1_c3} }\label{spp_1}
\end{align}
where the nonanticipatory constraint \eqref{p1_c4} is removed. $ Q_{t-1}(e_{t-1})$ is the maximized operating profit given the beginning state $e_{t-1}$, and $V_t(\cdot)$ is the value function defined as the expectation of the maximized arbitrage profit over the next time period
\begin{align}
    V_{t}(e_{t}) &= 
    \mathbb{E}\Big[ Q_{t}(e_{t}) \Big] \label{spp_vf}\,.
\end{align}
It is clear that $V_t(e_t)$ is valued and differentiable over the SoC range $[0,\;E]$. In addition, we define $v_t$ as the derivative of $V_t$, which indicates the SoC marginal value for a given SoC level. $v_t$ is mathematically expressed as
\begin{align}
    v_{t}(e_{t}) &= \begin{cases}
    -\infty & \text{if $e_t > E$} \\
    \mathbb{E}\Big[ q_{t}(e_{t}) \Big] & \text{if $0 \leq e_t \leq E$} \\
    \infty & \text{if $e_t < 0$}
    \end{cases}
\end{align}
where $q_t(\cdot) = \dot{Q}(\cdot)$ is the derivative of $Q(\cdot)$, hence the middle case ($0\leq e_t \leq E$) is obtained by moving the derivative operation inside the expectation in \eqref{spp_vf}. For the ease of presentation, we extend the range of $v_{t}(e_{t})$ beyond $[0,\;E]$ and assign infinity values to $v_t$ to model the infeasibility when $e_t > E$ or $e_t<0$. This representation also ensures a momentum to converge the SoC inside the feasible range. The impact of this infinity value definition in theoretical derivation and practical implementation will be further discussed in this paper.
\end{subequations}

\section{Main Results}
We derive an analytical formulation for calculating the SoC marginal value at a particular energy level $v_t(e)$  from the price distribution function and the value function for the next time period. Hence the entire SoC marginal value function can be obtained by sampling SoC over the feasible range $[0,E]$, and the result can be used for calculating the current value function via piece-wise linear approximations. Thus we obtain the storage value at the current time interval and also its optimal control result once a new price signal is received. The analytical calculation has excellent computational speed allowing detailed sampling of SoC, achieving high accuracy under the use of piece-wise linear approximations.

We start with a nontraditional approach by directly quantifying the cumulative distribution of $q_t(e)$, i.e., the marginal maximized arbitrage profit with respect to SoC at a given starting SoC $e$ over time period $t$. The SoC marginal value function $v_t$ is then calculated as the expectation of $q_t(e)$ using the obtained distribution following the standard expectation calculation as
\begin{align}
    v_t(e) = \int_{-\infty}^{\infty} x \mathbf{Pr}[q_t(e) = x] dx
\end{align}
where $\mathbf{Pr}[q_t(e) = x]$ is the distribution function of $q_t(e)$. At the end of this section, we introduce a piece-wise linear numerical algorithm for the recursive calculation of $v_t$. 

\subsection{Dual Decomposition}
We first convexify the non-simultaneous charging and discharging constraint \eqref{p1_c5} into
\begin{align}
   p\up{d} \geq 0 \text{ only if $\lambda_t > 0$} \label{p1_c6}
\end{align}
based on the observation that a sufficient condition for simultaneous charging and discharging is due to the occurrence of negative prices~\cite{xu2017scalable} and the intuitive observation that energy storage should never discharge when the price is negative. This setting is also enforced in most market designs that require selling offers to be non-negative. We will also discuss the validity of this convex relaxation later in Remark~\ref{re:nc} from a dynamic programming perspective. 
\begin{remark}\label{remark_concave}(\emph{Convex optimization})
 After the convex relaxation in \eqref{p1_c6}, the objective \eqref{spp_1} maximizes a concave function where $V_t$ is concave, and its derivative $v_t$ is a monotonic decreasing function.
\end{remark}

We apply dual decomposition~\cite{boyd2007notes} to the stage-wise decomposed problem in \eqref{spp_1} using the following dual variable $\theta_t$ associated with the SoC evolution constraint \eqref{p1_c1}, which we restate below:
\begin{subequations}
\begin{align}
    \theta_{t} : e_t - e_{t-1} = -p\up{d}_t/\eta + p\up{c}_t\eta 
\end{align}
and decomposes \eqref{spp_1} into 
\begin{align}
   \max_{p\up{c}_t, p\up{d}_t} \;&\lambda_t (p\up{d}_t-p\up{c}_t) - cp\up{d}_t - \frac{\theta_t}{\eta} p\up{d}_t + \theta_t\eta p\up{c}_t \text{ s.t. \eqref{p1_c2}, \eqref{p1_c6}}\label{dd_p1} \\
    \max_{e_t} \;&V_t(e_t) - \theta_te_t \text{ s.t. \eqref{p1_c3}}\label{dd_p2}
\end{align}
with the dual updating rule with step size $\epsilon\in\mathbb{R}^+$ 
\begin{align}
    ({1}/{\epsilon})\dot{\theta}_{t-1} = e_t - e_{t-1} + {p\up{d}_{t}}/{\eta} - p\up{c}_{t} \eta\label{dd_mc}
\end{align}
which simplifies the decision space and the treatment of the binary charge and discharge logic in \eqref{p1_c5}. The update rule \eqref{dd_mc} will be used in our later proofs, but will not be used in the proposed algorithm.
\end{subequations}
\begin{proposition}(\emph{Dual decomposition})\label{pro_dd}
\begin{subequations}
The solution to the dual decomposed problem \eqref{dd_p1} is
\begin{align}
    p\up{d}_t &= \begin{cases}
    P & \text{if $\lambda_{t} > [\theta_t/\eta + c]^+$} \\
    0 & \text{else}
    \end{cases} \label{pro1_e1}\\
        p\up{c}_t &= \begin{cases}
    P & \text{if $\lambda_{t} < \theta_t\eta$} \\
    0 & \text{else}
    \end{cases} \label{pro1_e2}
\end{align}
and for \eqref{dd_p2}
\begin{align}
    e_t = \begin{cases}
    E & \text{if $\theta_t < v_t(E)$} \\
    v^{-1}_t(\theta_t) & \text{if $v_t(E) \leq \theta_t \leq v_t(0)$} \\
    0 & \text{if $\theta_t > v_t(0)$} \\
    \end{cases} \label{pro1_e3}
\end{align}
where $ v^{-1}$ is the inverse function of $v_t$, and $[x]^+ = \max\{x,0\}$ is the positive value function.
\end{subequations}
\end{proposition}
The result in this proposition is obtained utilizing the first order optimality condition for convex optimization and limiting the result inside the upper and lower bound constraints. The condition $\lambda_{t} > [\theta_t/\eta + c]^+$ for $p\up{d}_t = P$ is equivalent to enforcing $\lambda_{t} > \theta_t/\eta + c$ (first-order optimality condition) and $\lambda_t > 0$ (constraint \eqref{p1_c6}). Note that although the actual solution $p\up{d}_t$ and $p\up{c}_t$ to \eqref{spp_1} may take on any value between the range $[0,P]$ since $V_t$ is a continuous or piece-wise linear function. But after the dual decomposition, $p\up{d}_t$ and $p\up{c}_t$ will either be $0$ or $P$, as \eqref{dd_p1} is a linear problem (recall $\theta_t$ is treated as a constant here) thus the solution must fall into one of the feasible region vertexes, based on the simplex algorithm. This is also an important characteristic of the dual decomposition that the solution may not be feasible, i.e., even if we plug in the optimal value of $\theta_t$, the decomposed result may not satisfy the master constraint in \eqref{dd_mc}.

In addition, the dual decomposition offers an insight of understanding simultaneous charging and discharging from a dynamic programming perspective:
\begin{remark}\label{re:nc}\emph{(Necessary conditions for simultaneous charging and discharging)}
We can  tell from \eqref{pro1_e1} and \eqref{pro1_e2} that in the cases of simultaneous charging and discharging (if we relax constraint \eqref{p1_c5} ), it must follow that $\lambda_t/\eta + c < \theta_t < \lambda\eta$. Since $c>0$, then a necessary condition for this inequality relationship is $\theta_t < 0$ and $\lambda_t < 0$. Hence if we enforce $p\up{d}_t = 0$ when $\lambda_t < 0$, we can avoid simultaneous charging and discharging, which also follows a straight forward intuition that the storage should not sell energy when price is negative. Also note that if $\theta_t \geq 0$, then simultaneous charging and discharging will not occur.
\end{remark}

\subsection{Quantifying SoC Value Distribution}
Based on the dual decomposition result, we obtain the main theorem  that quantifies the cumulative distribution of $q_t$ which is the marginal value of stored energy (SoC):
\begin{theorem}\label{th:soc}(\emph{SoC marginal value distribution})
The cumulative distribution of $q_{t-1}(e)$ can be quantified from the value function derivative of the next time period $v_t$ as ($x\in\mathbb{R}$)
\begin{align}
        \mathbf{Pr}[&q_{t-1}(e) \leq x] = \nonumber\\
        &\begin{cases}
        0 \quad & \text{if $x < v_t(e +P\eta)$} \\
        F_t({x}{\eta}) & \text{if $ v_t(e+P\eta) \leq x < v_t(e)$} \\
        F_t([{x}/{\eta}+c]^+) & \text{if $ v_t(e) \leq x < v_t(e-P/\eta)$} \\
        1 & \text{if $x \geq v_t(e-P/\eta)$}
        \end{cases} \label{th1}
\end{align}
\end{theorem}
The key idea of this proof is utilizing the analytical dual decomposition result in Proposition~\ref{pro_dd}, with which the dual value $\theta_t$ can be updated based on \eqref{dd_mc} and an inequality relationship can be established between the dual variable and the ending SoC of the considered time period. The detailed proof of this theorem is deferred to Appendix~A.

\subsection{Value Function for Risk-Neutral Policies}
Based on Theorem~\ref{th:soc}, the expected SoC marginal value function $v_t$ is calculated in the following proposition:
\begin{proposition}\label{pro:v}(\emph{Value function derivative})
The risk-neutral value function derivative can be calculated from the distribution function $f_t$ and $F_t$, the value function $v_t$, power rating $P$, and efficiency $\eta$ as 
\begin{align}
    v_{t-1}&(e) = \mathbb{E}[q_{t-1}(e)] =\nonumber\\
    & v_t(e+P\eta)F_t\big(v_t(e+P\eta)\eta\big) \nonumber\\
    & + \frac{1}{\eta}\int_{v_t(e+P\eta)\eta}^{v_t(e)\eta} uf_t(u)\;du\nonumber\\
    & + v_t(e)\Big[F_t\big([v_t(e)/\eta + c]^+\big) -  F_t\big(v_t(e)\eta\big)\Big] \nonumber\\
    & + \eta\int_{[v_t(e)/\eta  +c]^+}^{[v_t(e-P/\eta)/\eta + c]^+} w f_t(w)\;dw \nonumber\\
    & - c\eta\Big[F_t\big({[v_t(e-P/\eta)}/{\eta} + c]^+\big) - F_t\big([{v_t(e)}/{\eta} + c]^+\big)\Big]\nonumber\\
    & + v_t(e-P/\eta)\Big[1-F_t\big([v_t(e-P/\eta)/\eta + c]^+\big)\Big]\label{pro1}\,.
\end{align}
\end{proposition}

% This result can be written out with respect to the distribution used to model the forecast error. In the example of using normal distribution, $v_{t-1}(e)$ can then be written as
% \begin{corollary}
% If $\lambda_t$ follows normal distribution, then $v_{t-1}(e)$ can be calculated as
% \begin{align}
%     v_{t-1}(e) =&  - \frac{1}{\eta}\sigma_t^2 \Big[f_t\big(v_t(e)\eta\big)-f_t\big(v_t(e+P\eta)\eta\big)\Big] \nonumber\\
%     &-\eta \sigma_t^2 \Big[f_t\big(v_t(e-P/\eta)/\eta\big)-f_t\big(v_t(e)/\eta\big)\Big]
% \end{align}
% \end{corollary}
Proof of this proposition is deferred to Appendix~B. \eqref{pro1} can be calculated in closed form except for the second and the fourth term which involves integration, however these two terms are identical to calculating conditional expectations as (let $a\leq b$, $(a,b)\in\mathbb{R}^2$)
\begin{align}
\int_{a}^bxf(x)\;dx = {E[x|a\leq x \leq b]}\big(F(b)-F(a)\big)   
\end{align}
which can be evaluated efficiently for most known distributions. Alternatively, the integral calculation can be discretized into summation if the price distribution is described via discrete samples. Another difficulty in evaluating \eqref{pro1} recursively is calculating $v_t$,
for this we demonstrate a sampling algorithm in the next section.

\subsection{Algorithmic Implementation}
We discretize $v_t$ by modeling it as an vector $\{v_{t,j}\}$ in which each element is associated with equally spaced SoC samples $\{e_j = (j-1)\Delta e | j\in \{1,\dotsc, J\}\}$, where $\Delta e$ is the sample granularity, and the number of samples  $J = 1 + E/\Delta e$. For any value inquiry $v_t(e)$, we round $e$ to the nearest SoC samples and return the corresponding value. This discrete value function derivative $\tilde{v}_t$ is formally defined as
\begin{align}
    \tilde{v}_t(e) = v_{t,j},\;  j = \mathbf{proj}_j(e) 
\end{align}
where the index projection function $\mathbf{proj}_j(e)$ rounds the index of the SoC samples $e_j$ to the closest to $e$. An algorithmic evaluation of $v_t$ as in \eqref{pro1} can thus be achieved using $\tilde{v}_t$.
\begin{remark}(\emph{Piece-wise linear value function})
Note that  discretizing the value function derivative  $v_t$ is equivalent as to approximating $V_t$ using piece-wise linear functions since $V_t$ is the integral of $v_t$. 
\end{remark}
\begin{remark}(\emph{Complexity analysis})
The proposed algorithm achieves linear time complexity and constant space complexity with respect to the look-ahead horizon $T$. If we assume $J$ number of SoC samples, then at each time interval, we need to execute \eqref{pro1} for $J$ times. Thus over a time horizon of $T$, \eqref{pro1} will be executed a total of $TJ$ times, hence a linear time complexity. At each time interval, we only need to record $v_t$ for the calculation of $v_{t-1}$, thus the calculation space needed will not increase with $T$, hence constant space complexity.
\end{remark}
\subsection{Extension to General Objectives}
We discuss how Theorem~\ref{th:soc} and Proposition~\ref{pro:v} could be extended to solve stochastic storage control for maximizing a general concave objective function, in which the objective function in \eqref{p1_obj} is rewritten as
\begin{align}
    \max_{p_t} \mathbb{E}\Big[\sum_{t=1}^T R_t(p\up{d}_t - p\up{c}_t) \Big] + V_T(e_T)
\end{align}
where $R_t(x)$ is a scalar concave function that models the revenue received from the market via the action $x$. The main difficulty in dealing with a general objective function is that we can no longer discretize the action space precisely according to the dual decomposition result in \eqref{dd_p1} and \eqref{dd_p2}. For example, if $R_t(x) = a x^2 + bx$ , then the solution to $p\up{d}_t$ after the dual decomposition according to first order optimality condition is
\begin{align}
    p\up{d}_t = \frac{\theta_t - b\eta}{2a\eta}
\end{align}
which is a linear function with respect to $\theta_t$.

In this case, Theorem~\ref{th:soc} and Proposition~\ref{pro:v} can be applied if we discretize the action space $(p\up{d}_t, p\up{c}_t)$ after taking the dual decomposition. This is equivalent to approximating $R_t$ using piece-wise linear functions. For example, if we assume $R_t$ is associated with slopes $c_j$ and segments $[P_{j-1}, P_j]$ as
\begin{align}
    \dot{R}_t(x) = c_j \quad \text{if $P_{j-1} \leq x < P_j$}
\end{align}
then Theorem~\ref{th:soc}  can be expanded to the following format
\begin{align}
    &\mathbf{Pr}[q_{t-1}(e) \leq x] = F_{t,j}([x/\eta + c_j]^+) \\
    &\quad \text{if $v_t(e-P_{j-1}/\eta) \leq x \leq v_t(e-P_j/\eta)$}\nonumber
\end{align}
where $F_{t,j}$ is the distribution of the $j$th cost segment $c_j$, and Proposition~\ref{pro:v} can be applied similarly. The difficulty in this approach is designing a high dimensional forecast for the cost segment $c_j$ which is beyond the scope of our discussion in this paper.

\section{Empirical Analysis}
% We first demonstrate examples of energy storage valuation under different price scenarios and discuss the correlation of the storage value with the storage size and the price standard deviation. All batteries are assumed to have 90\% round trip efficiency. Second, we also include an electric vehicle charging example using historical price data from New York City\footnote{Data available at the New York ISO dataport: \url{https://www.nyiso.com/}.}. 

\subsection{Computation Performance}
The proposed algorithm is implemented in Matlab~\cite{matlab} and is compared to solving problem \eqref{p1} using the FAST stochastic dual dynamic programming (SDDP) toolbox in Matlab~\cite{fast} with Gurobi~\cite{gurobi} embedded, both were performed on a MacBook Pro with 2.3~GHz CPU and 16~GB memory. The SDDP solver is set to simulate 25 Monte-Carlo scenarios at each forward pass and the maximum iteration number is 20.

% \begin{figure}[ht]%
% 	\centering
% 	\subfloat[Comparison of control results in a random scenario.]{
% 		\includegraphics[trim = 10mm 0mm 10mm 0mm, clip, width = .95\columnwidth]{sddp_power.eps}
% 		\label{Fig:size1}%
% 	}
% 	\\
% 	\subfloat[SDDP convergence result (minimize the negative profit).]{
% 		\includegraphics[trim = 10mm 0mm 10mm 0mm, clip, width = .95\columnwidth]{sddp_iter.eps}
% 		\label{Fig:size2}%
% 	}
%     \caption{Comparison with SDDP with 10 stages and each stage has 10 scenarios. The SDDP solver is able to converge in 14 iterations, with a total computation speed of 57 seconds. The proposed algorithm finished in 23~ms.}%
%     \label{Fig:size}
% \end{figure}

\begin{figure}[t]%
	\centering
	\subfloat[Solution time of the SDDP solver in minutes (20 iterations max).]{
		\includegraphics[trim = 10mm 0mm 10mm 0mm, clip, width = .95\columnwidth]{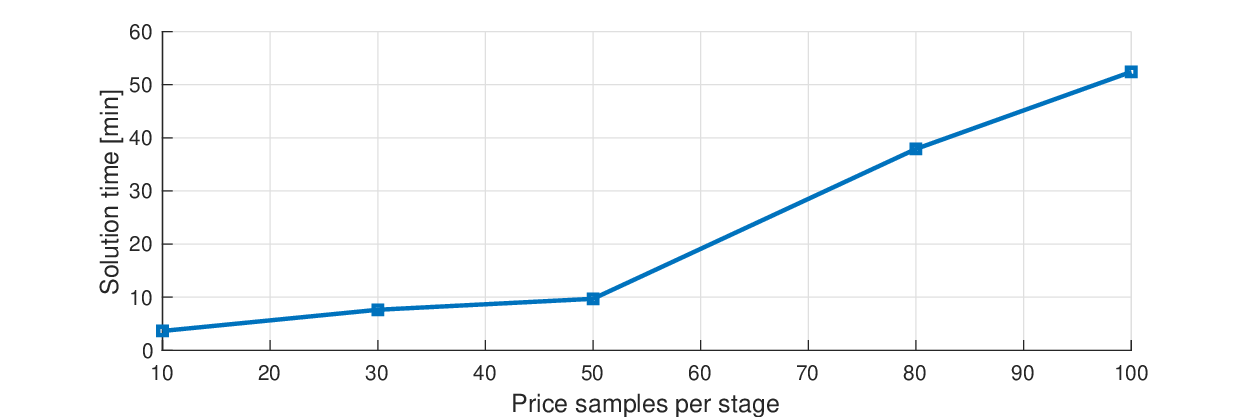}
		\label{Fig:comp1}%
	}
	\\
	\subfloat[Solution time of the proposed method in milliseconds.]{
		\includegraphics[trim = 10mm 0mm 10mm 0mm, clip, width = .95\columnwidth]{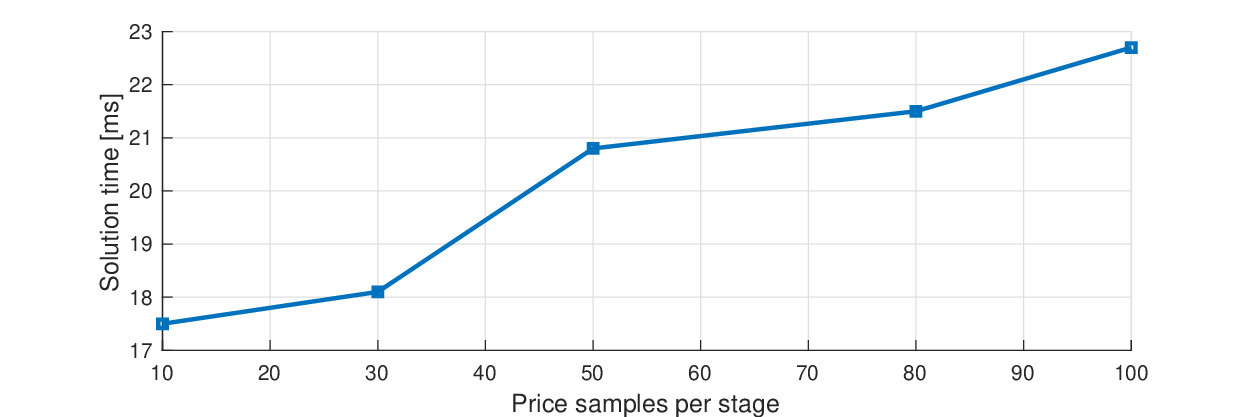}
		\label{Fig:comp3}%
	}
	\\
	\subfloat[Average objective value in Monte-Carlo simulations.]{
		\includegraphics[trim = 10mm 0mm 10mm 0mm, clip, width = .95\columnwidth]{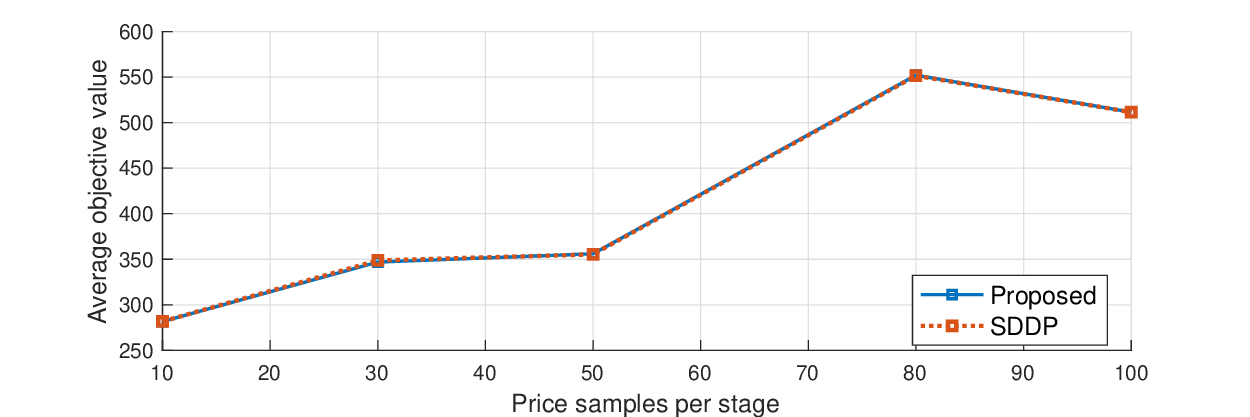}
		\label{Fig:comp2}%
	}
    \caption{Comparison between the proposed method and the benchmark SDDP solver over 24 look-ahead stages.  }%
    \label{Fig:comp}
\end{figure}

We consider 24 uncertainty stages and test both methods considering different number of distribution samples (nodes) per stage, the result in shown in Fig.~\ref{Fig:comp}. For example, if we consider 10 samples per stage, the will be a total of $10^{24}$ scenarios. The computation time of SDDP increases significantly with the number of samples considered, while the proposed method finishes within 25 milliseconds in all test cases (shown in Fig.~\ref{Fig:comp3}), which is up to 100,000 times faster than SDDP in the simulated cases. The optimality of the solution is tested using Monte-Carlo simulation by sampling different scenarios. As shown in Fig.~\ref{Fig:comp2}, both methods provides very similar results on the average objective function value.

% Both methods produced similar results in small test cases, i.e., 10 stages each with 10 scenarios (a total of $10^9$ scenarios as the first stage is deterministic). The FAST toolbox is set to perform 25 Monte-Carlo scenario samples at each forward path. In these small test cases, the proposed method finished within 20 milliseconds, while the FAST toolbox can converge within a few minutes using less than 20 iterations. Objective function value is tested using Monte-Carlo simulation, the error between the two method is within 1\%, with the proposed method offering less total profit due to the error introduced in the piece-wise linear approximation. 

% However, when increasing the test size, either with a longer forecast horizon or with more scenarios at each stage, the solution time of SDDP will increase significantly.

\subsection{Impact of Price Standard Deviations}

\begin{figure}[ht]%
	\centering
	\subfloat[Example day-ahead price and forecast distributions.]{
		\includegraphics[trim = 10mm 0mm 10mm 0mm, clip, width = .95\columnwidth]{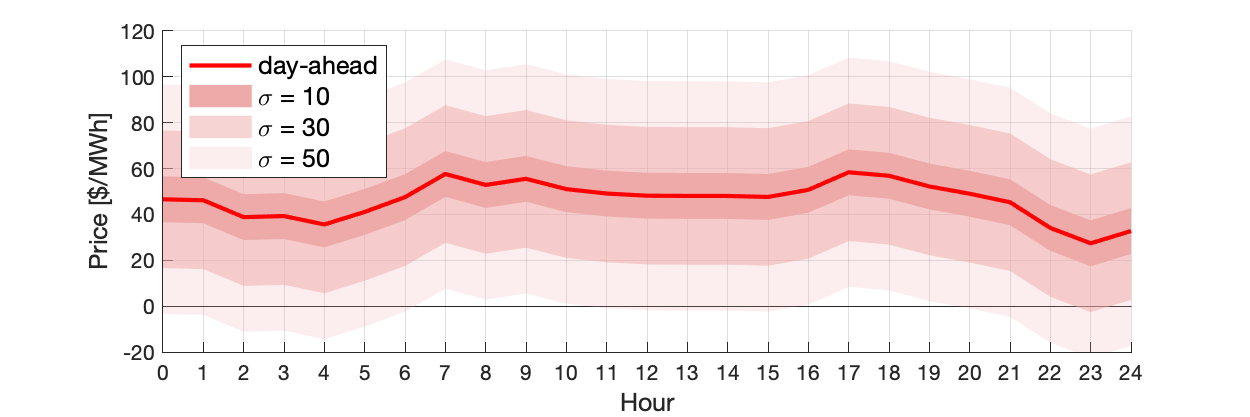}
		\label{Fig:dev1}%
	}
	\\
	\subfloat[Energy storage marginal value ranges.]{
		\includegraphics[trim = 10mm 0mm 10mm 0mm, clip, width = .95\columnwidth]{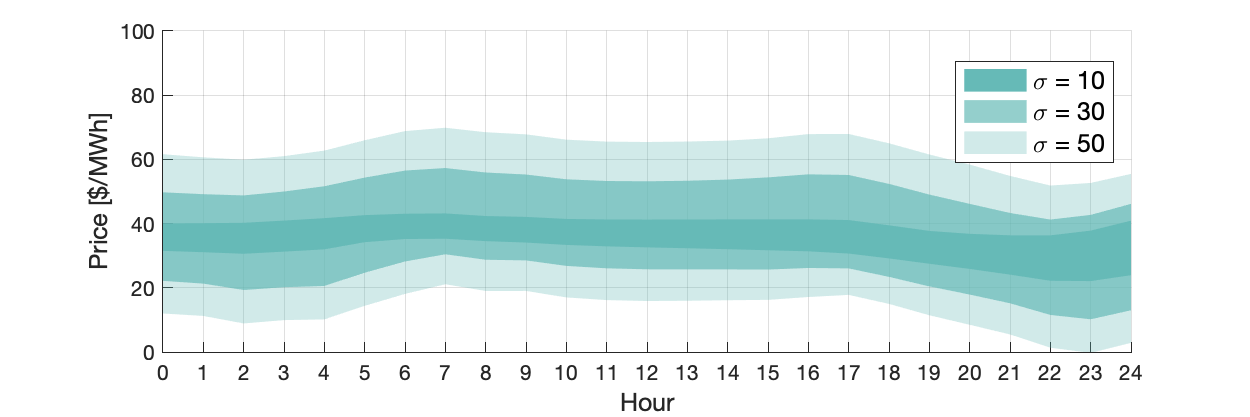}
		\label{Fig:dev2}%
	}
    \\
	\subfloat[Storage value for different SOCs at hour 12.]{
		\includegraphics[trim = 10mm 0mm 10mm 0mm, clip, width = .95\columnwidth]{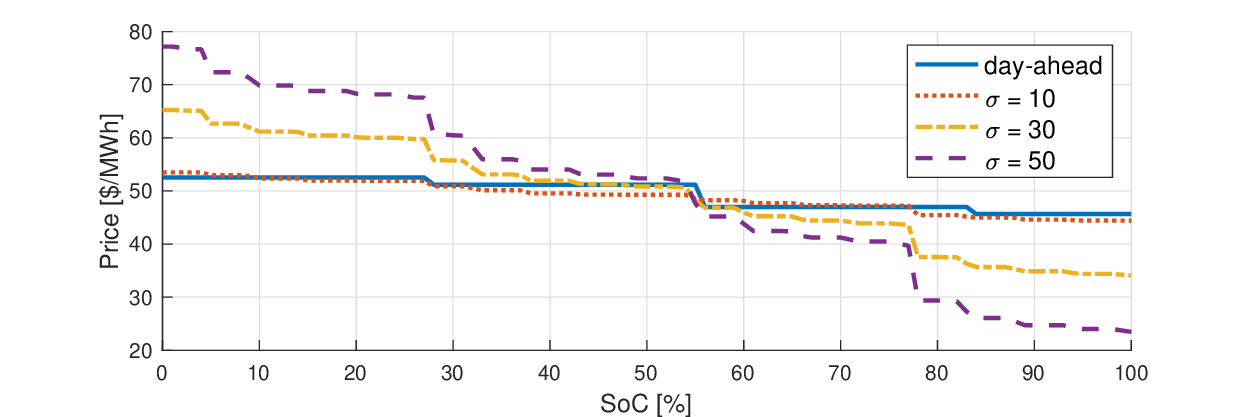}
		\label{Fig:dev3}%
	}
    \caption{Example of energy storage value vs. different price distributions, assuming normal distribution with standard deviations $\sigma$ of 10, 30, and 50. }%
    \label{Fig:dev}
\end{figure}
We demonstrate the impact of price uncertainty in the storage valuation. We model the real-time price forecast as zero-mean normal distributions imposed over the day-ahead price results, and include three cases with the standard deviation $\sigma$ of 10, 30, and 50, as shown in Figure~\ref{Fig:dev1}.

We value a 4-hour energy storage device (a fully charged storage will take 4 hours to fully discharge at rated power) using the forecast information.
Figure~\ref{Fig:dev2} shows the the resulting value price range (i.e., $v_t$), the upper edge indicates the marginal storage value at 0\% SoC, and the lower edge at 100\% SoC. The valuation also considers operation beyond the 24 hour forecast horizon so that the storage value does not go to zero at the end of the operation. As the price uncertainty increases, the storage value  also spans a wider range. The results are also illustrated in Figure~\ref{Fig:dev3}, in which the storage values are plotted versus different SoC. 

\subsection{Storage Charging Case Study}

\begin{figure}[ht]%
	\centering
	\subfloat[Day-ahead price, real-time price, price forecast uncertainty (boxes in the plot), New York City, Feb.1 2018.]{
		\includegraphics[trim = 10mm 0mm 10mm 0mm, clip, width = .95\columnwidth]{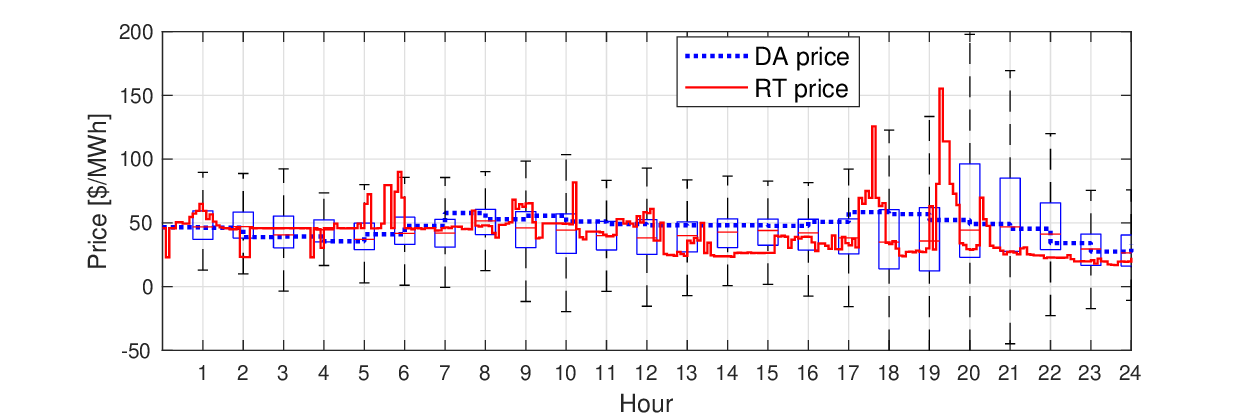}
		\label{Fig:ny1}%
	}
	\\
	\subfloat[SoC value ranges using different price forecasts. ]{
		\includegraphics[trim = 10mm 0mm 10mm 0mm, clip, width = .95\columnwidth]{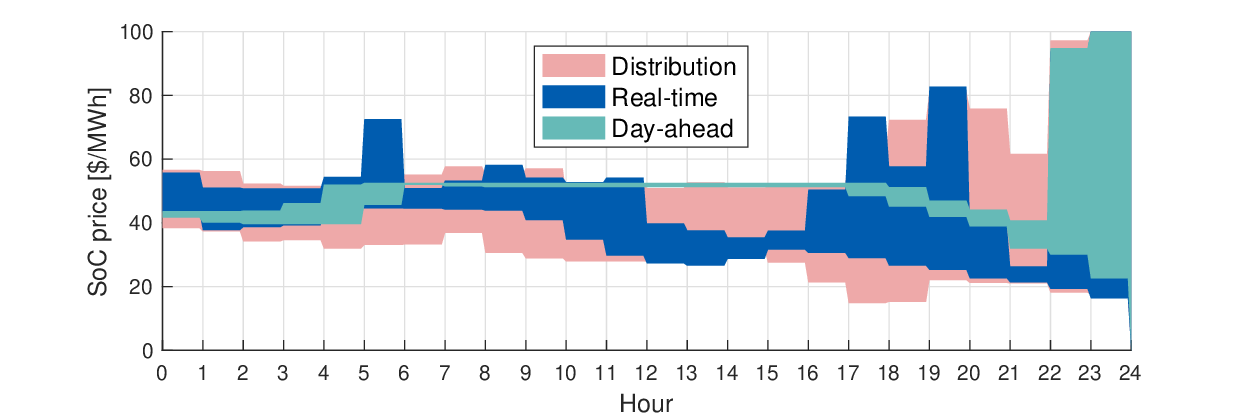}
		\label{Fig:ny2}%
	}
	\\
	\subfloat[SOC pattern based on different price forecasts.]{
		\includegraphics[trim = 10mm 0mm 10mm 0mm, clip, width = .95\columnwidth]{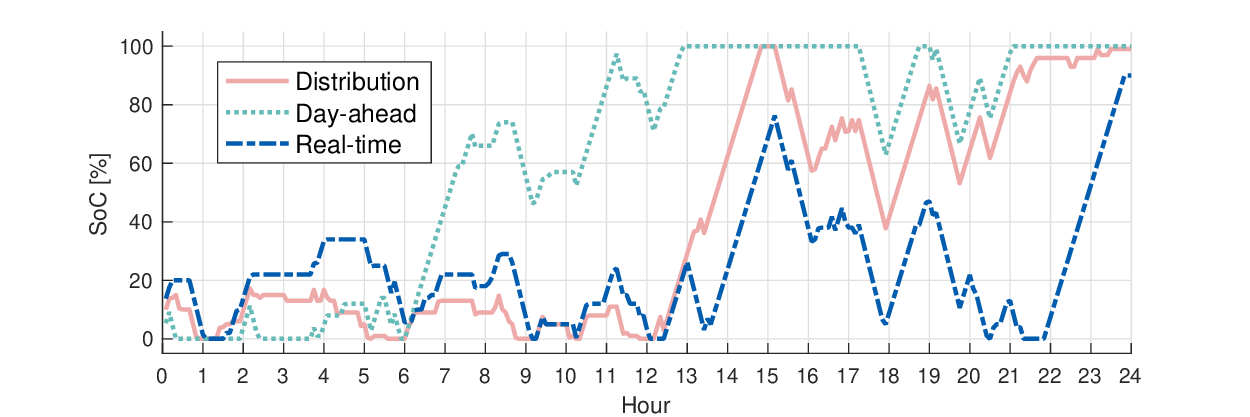}
		\label{Fig:ny3}%
	}
    \caption{Charging a 100kW/200kWh battery from 10\% to 90\% SoC using real-time prices and different price forecasts. }%
    \label{Fig:ny}
\end{figure}

We consider an example application in which we want to charge a 100kW/200kWh storage device from 10\% SoC to 90\% SoC via real-time prices under the least cost or even earning a profit (negative cost). This is modeled into \eqref{p1} by setting $e_0$ to 10\%, and derivative of final value  function $v_T$ is modeled as a step function with a value of 100\$/MWh from 0\% SoC to 90\% SoC, and zero value after 90\% SoC. This application closely resembles the charging of an electric vehicle if exposed to real-time price uncertainties. We use the New York City price data as shown in Figure~\ref{Fig:ny1}\footnote{Data available at the New York ISO dataport: \url{https://www.nyiso.com/}.}, including the day-ahead (DA) price, and real-time (RT) price, and the probability forecast of the RT price based on historical DA-RT price biases over January 2018. Using the proposed algorithm, we calculate the hourly SoC value curves using the DA price, real-time price, and the price distribution, respectively. The resulting valuation ranges are shown in Figure~\ref{Fig:ny2}. Note that the valuation using the RT price is a perfect forecast benchmark as it is impossible to know real-time price realizations beforehand, while the valuation using DA prices is a simple heuristic assuming the RT price will be the same as the DA price.

We then perform the three different storage control strategies using the RT price series and the calculated value functions. The resulting SoC series are shown in Figure~\ref{Fig:ny3}. In all three cases, the storage reaches the target 90\% SoC level in the end, while also gaining some profit. The total market profit using DA price-based valuation is \$2 and for RT price-based valuation is \$8. Our proposed price distribution-based valuation results in a revenue of \$4, which is substantially better than the DA case but is of course less than the perfect information RT case.

\section{Conclusion}

We presented a framework for valuation of energy storage operation by solving a multi-stage energy arbitrage problem under price uncertainty. Our proposed solution has very high computational speed as it involves only basic computational operations. The computation speed of the proposed method significantly surpassed the benchmark SDDP solver and it requires only basic arithmetic and logic operations. Notably, the SDDP solver will be very inefficient to solve storage operation over high time resolutions (i.e., real-time market arbitrage at 5-minute granularity). These cases will have hundreds of uncertainty stages and the SDDP solver will take several hours to finish due to its super-linear time complexity, while the proposed method is still guaranteed to finish within a few seconds.  

The framework can help energy storage participants in calculating their storage value instantaneously based on the most recent price forecast. Moreover, our algorithm can be implemented using very simple hardware and software which is ideal for optimizing distributed energy storage devices. Our future work includes expanding the proposed method to coordinating the control of multiple storage devices such as EV charging, and research the impact of inter-temporal price correlation to the valuation.

\bibliographystyle{IEEEtran}	% (uses file "plain.bst")
\bibliography{IEEEabrv,literature}		% expects file "myrefs.bib"

\appendix

\subsection{Proof of Theorem 4}
% \begin{proof}
% We first apply dual decomposition~\cite{boyd2007notes} to \eqref{spp_1} by making the dual variable $\theta_t$ as a master problem variable ( treated as a constant in the subproblems) and apply Lagrangian relaxation to \eqref{p1_c1} and write out two subproblems with respect to $p_t$ and $e_t$ separately:
\begin{subequations}
% \begin{align}
%     \max_{p_t\in [-P,P]} &\lambda_t p_t - c[p_t]^+ - \theta_t([p\up{\pi}_t]^+/\eta + [p\up{\pi}_t]^-\eta) \\
%     \max_{e_t\in [0,E]} &V_t(e_t) - \theta_te_t
% \end{align}
% in which the subproblem solution can be expressed as functions of $\theta_t$ as 
% \begin{align}
%     p^*_t(\theta_t) &= \begin{cases}
%     P \; & \text{if $\lambda_t \geq \theta_t/\eta + c$} \\
%     0 \; & \text{if $\theta_t\eta < \lambda_t < \theta_t/\eta + c$} \\
%     -P \; & \text{if $\lambda_t \leq  \theta_t\eta$} 
%     \end{cases}\label{th1_a}\\
%     e_t^*(\theta_t) &= [v_t^{-1}(\theta_t)]^{E}_0\label{th1_b}
% \end{align}
% and the dual variable $\theta_t$ can be iteratively solved with the following master updating algorithm 
% \begin{align}
%     \frac{\dot{\theta}_{t-1}}{\epsilon} = e^*_t(\theta_t) - e_{t-1} + \frac{[p^*_t(\theta_t)]^+}{\eta} + [p^*_t(\theta_t)]^-\eta
% \end{align}
% where $\epsilon$ is a non-negative step size which may be different at each iteration $k$. 
For the convenience of presentation, we denote $p\up{d}_t(x)$, $p\up{c}_t(x)$, $e_t(x)$ as the solution to the dual decomposed problem when the dual variable value is set to $x$, according to Proposition~\ref{pro_dd}. 
Then for all $x$ greater than the optimal dual variable value, we have $x \geq \theta_t$, $x\in\mathbb{R}$ 
\begin{align}
    e_t(x)  \leq e_{t-1} + \frac{p\up{d}_t(x)}{\eta} - p\up{c}_t(x)\eta
\end{align}
For simplicity, we denote the right-hand side of the above equation as $e\up{p}_t(x)$, which is the ending SoC resulting from the power sub-problem ~\eqref{dd_p1}. Now according to the dual decomposition master constraint updating rule in \eqref{dd_mc}, we reformulate the previous equation and apply $v_t(\cdot)$ to both sides 
\begin{align}
    e_t(x) &\leq e\up{p}_t(x) \\
    v_t(e_t(x)) &\geq v_t(e\up{p}_t(x) \\
    x &\geq v_t (e\up{p}_t(x)) \label{th1_c}
\end{align}
according to \eqref{pro1_e3} from Proposition~\ref{pro_dd} and Remark~\ref{remark_concave} that $v_t(\cdot)$ is a non-increasing function. Meanwhile, $e\up{p}_t(x)$ can be analytically written-out by substituting \eqref{pro1_e1} and \eqref{pro1_e2} into:
\begin{align}
    e\up{p}_t(x) = \begin{cases}
        e_{t-1} - P/\eta & \text{if $ \lambda_t > [x/\eta + c]^+$} \\
        e_{t-1} & \text{if $x\eta \leq \lambda_t \leq [x/\eta + c]^+$} \\
        e_{t-1} + P\eta & \text{if $\lambda_t < x\eta$}
    \end{cases}\,. \label{th1_d}
\end{align}

Now given any value $x$, we can determine its inequality relationship to $\theta_t$ by combining \eqref{th1_c} and \eqref{th1_d} as
\begin{align}
    \mathbf{Pr}[&q_{t-1}(e_{t-1}) \leq x] =\nonumber\\
    &\mathbf{Pr}[\lambda_t > [x/\eta + c]^+ ]\cdot \mathbf{1}[v_t(e_{t-1} - P/\eta) \leq x] \nonumber\\
    \quad + &\mathbf{Pr}[x\eta \leq \lambda_t [\leq x/\eta + c]^+]\cdot \mathbf{1}[v_t(e_{t-1}) \leq x] \nonumber\\
    \quad +  &\mathbf{Pr}[\lambda_t < x\eta]\cdot \mathbf{1}[v_t(e_{t-1}+ P\eta) \leq x]
\end{align}
% where the first probability can be separated as
% \begin{align}
%     &\mathbf{Pr}[\lambda_t > x/\eta + c \textbf{ and } \lambda_t > 0] = \nonumber\\
%     & \mathbf{Pr}[\lambda_t > x/\eta + c]
% \end{align}
where the indicator function $\mathbf{1}[x]= 1$ if statement $x$ is true and $\mathbf{1}[x]= 0$ otherwise. This theorem is thus proved by substituting $f_t$ into as:
\begin{align}
    % \mathbf{Pr}[\lambda_t > 0] &= 1-F_t(0) \\
    % \mathbf{Pr}[\lambda_t > x/\eta + c ] &= 1 - F_t(x/\eta + c) \\
    \mathbf{Pr}[\lambda_t > [x/\eta + c]^+ ] &= 1-F_t([x/\eta + c]^+)\nonumber\\
    \mathbf{Pr}[x\eta \leq \lambda_t \leq [x/\eta + c]^+] & = F_t([x/\eta + c]^+) - F_t(x\eta)\nonumber\\
    \mathbf{Pr}[\lambda_t < x\eta] &= F_t(x\eta)\nonumber
\end{align}
where depending on the range of $x$, we can combine terms within the same range which gives us the result in \eqref{th1}.
The subscript of $e_{t-1}$ can thus be removed as it is the only SoC variable in this equation.  
\end{subequations}
% \end{proof} 

\subsection{Proof of Proposition 2}

\begin{proof}
We start with the standard expectation calculation which is written in the following form
\begin{subequations}
\begin{align}
    v_{t-1}(e) = \mathbb{E}\big[ q_{t-1}(e) \big] 
    & =  \int_{-\infty}^{\infty} x\frac{\partial \mathbf{Pr}[q_{t-1}(e) \leq x]}{\partial x} \; dx \nonumber\\
    & = \int_{-\infty}^{\infty} x\mathbf{Pr}[q_{t-1}(e) = x]\;dx\label{pro1_1}
    % &=  \int_{v_t(e+P\eta)}^{v_t(e)} {x}{\eta}f_t({x}{\eta})\;dx   + \int_{v_t(e)}^{v_t(e-P/\eta)} \frac{x}{\eta} f_t(\frac{x}{\eta}+c)\;dx\nonumber\\
    % &=  \frac{1}{\eta}\int_{v_t(e+P\eta)\eta}^{v_t(e)\eta} xf_t(x)\;dx   \nonumber\\ &\quad + \eta\int_{v_t(e)/\eta}^{v_t(e-P/\eta)/\eta} x f_t(x)\;dx\nonumber\\
\end{align}
where based on Theorem~1 we write $\mathbf{Pr}[q_{t-1}(e) = x]$ as
\begin{align}
    &\mathbf{Pr}[q_{t-1}(e) = x] = \label{pro1_2}\\
        &\begin{cases}
        F_t(x\eta) & \text{if $x = v_t(e+P\eta)$}\\
        \eta f_t({x}{\eta}) & \text{if $ v_t(e+P\eta)< x < v_t(e)$} \\
        F_t([\frac{x}{\eta}+c]^+)-F_t(x\eta) & \text{if $x = v_t(e)$}\\
        \frac{1}{\eta}f_t([\frac{x}{\eta}+c]^+)\mathbf{1}[\frac{x}{\eta}+c > 0] & \text{if $ v_t(e) < x < v_t(e-P/\eta)$} \\
        1 - F_t(x/\eta + c) & \text{if $x = v_t(e-P/\eta)$} \\
        0 & \text{else}
        \end{cases} \nonumber 
\end{align}
in which the second, fourth, and the sixth cases are obtained by taking the derivative to cases in \eqref{th1} directly over the respective range. For the discontinuous points $x = v_t(e+P\eta)$, $x = v_t(e)$, $x = v_t(e-P/\eta)$, their probabilities are the difference between the two neighbouring cumulative distribution results. For example, at $x = v_t(e+P\eta)$ is connected to $\mathbf{Pr}[q_{t-1}(e) \leq x] = 0$ and $\mathbf{Pr}[q_{t-1}(e) \leq x] = F_t(x\eta)$, thus we have 
\begin{align}
    \mathbf{Pr}[q_{t-1}(e) = x] &= \mathbf{Pr}[q_{t-1}(e) \leq  x] - \mathbf{Pr}[q_{t-1}(e) < x]\nonumber\\
    & = F_t(x\eta)
\end{align}
To connect \eqref{pro1_2} to the result in \eqref{pro1}, it is trivial to see that the first, third, and fifth case of \eqref{pro1_2} correspond to the first, third, and sixth case in \eqref{pro1}, calculated by multiplying the $x$ value with the corresponding discrete probability. For the second case, we let $u = x\eta$ and use the substitution rule for definitive integrals which gives
\begin{align}
    \int_{v_t(e+P/\eta)}^{v_t(e)} x\eta f_t(x\eta)\;dx =  \frac{1}{\eta}\int_{v_t(e+P\eta)\eta}^{v_t(e)\eta} uf_t(u)\;du
\end{align}
hence the second term in \eqref{pro1}. For the fourth case, let $w = [x/\eta + c]^+$ and hence
\begin{align}
    & \int_{v_t(e)}^{v_t(e-P/\eta)} \frac{x}{\eta} f_t([\frac{x}{\eta}+c]^+)\mathbf{1}[\frac{x}{\eta}+c > 0]\;dx\nonumber\\
    =& \int_{[v_t(e)/\eta  +c]^+}^{[v_t(e-P/\eta)/\eta + c]^+} w f_t(w)\;dw \nonumber\\
  &  - c \int_{v_t(e)}^{v_t(e-P/\eta)} f_t(\frac{x}{\eta}+c)\mathbf{1}[\frac{x}{\eta}+c > 0]\;dx \nonumber\\
    = & \int_{[v_t(e)/\eta  +c]^+}^{[v_t(e-P/\eta)/\eta + c]^+} w f_t(w)\;dw \nonumber\\
   & - c\eta\Big[F_t\big([\frac{v_t(e-P/\eta)}{\eta} + c]^+\big) - F_t\big([\frac{v_t(e)}{\eta} + c]^+\big)\Big]
\end{align}
which equals the remaining (fourth and fifth term) in \eqref{pro1}.
\end{subequations}
\end{proof}

% that's all folks
\end{document}